\documentclass[reqno,a4paper]{amsart}
\usepackage{amssymb}
\usepackage[colorlinks=true,hypertex]{hyperref}
\allowdisplaybreaks[4]
\theoremstyle{plain}
\newtheorem{thm}{Theorem}
\theoremstyle{remark}
\newtheorem{rem}{Remark}
\DeclareMathOperator{\td}{d\mspace{-2mu}}
\numberwithin{equation}{section}
\date{Completed on Monday 21 July 2008}
\date{}

\begin{document}

\title[Bounds for the ratio of two gamma functions]
{Bounds for the ratio of two gamma functions---From Wendel's limit to Elezovi\'c-Giordano-Pe\v{c}ari\'c's theorem}

\author[F. Qi]{Feng Qi}
\address[F. Qi]{Research Institute of Mathematical Inequality Theory, Henan Polytechnic University, Jiaozuo City, Henan Province, 454010, China}
\email{\href{mailto: F. Qi <qifeng618@gmail.com>}{qifeng618@gmail.com}, \href{mailto: F. Qi <qifeng618@hotmail.com>}{qifeng618@hotmail.com}, \href{mailto: F. Qi <qifeng618@qq.com>}{qifeng618@qq.com}}
\urladdr{\url{http://qifeng618.spaces.live.com}}

\begin{abstract}
In the survey paper, along one of main lines of bounding the ratio of two gamma functions, we look back and analyse some known results, including Wendel's, Gurland's, Kazarinoff's, Gautschi's, Watson's, Chu's, Lazarevi\'c-Lupa\c{s}'s, Kershaw's and Elezovi\'c-Giordano-Pe\v{c}ari\'c's inequalities, claim, monotonic and convex properties. On the other hand, we introduce some related advances on the topic of bounding the ratio of two gamma functions in recent years.
\end{abstract}



\thanks{The first author was partially supported by the China Scholarship Council}

\thanks{This paper was typeset using \AmS-\LaTeX}

\maketitle

\tableofcontents

\section{Introduction}

Recall~\cite[Chapter~XIII]{mpf-1993} and~\cite[Chapter~IV]{widder} that a function $f$ is said to be completely monotonic on an interval $I$ if $f$ has derivatives of all orders on $I$ and
\begin{equation}
(-1)^{n}f^{(n)}(x)\ge0
\end{equation}
for $x \in I$ and $n \geq0$. The celebrated Bernstein-Widder Theorem~\cite[p.~160, Theorem~12a]{widder} states that a function $f$ is completely monotonic on $[0,\infty)$ if and only if
\begin{equation}\label{converge}
f(x)=\int_0^\infty e^{-xs}\td\mu(s),
\end{equation}
where $\mu$ is a nonnegative measure on $[0,\infty)$ such that the integral~\eqref{converge} converges for all $x>0$. This tells us that a completely monotonic function $f$ on $[0,\infty)$ is a Laplace transform of the measure $\mu$.
\par
It is well-known that the classical Euler's gamma function may be defined for $x>0$ by
\begin{equation}\label{egamma}
\Gamma(x)=\int^\infty_0t^{x-1} e^{-t}\td t.
\end{equation}
The logarithmic derivative of $\Gamma(x)$, denoted by $\psi(x)=\frac{\Gamma'(x)}{\Gamma(x)}$, is called the psi or digamma function, and $\psi^{(k)}(x)$ for $k\in \mathbb{N}$ are called the polygamma functions. It is common knowledge that the special functions $\Gamma(x)$, $\psi(x)$ and $\psi^{(k)}(x)$ for $k\in\mathbb{N}$ are fundamental and important and have much extensive applications in mathematical sciences.
\par
The history of bounding the ratio of two gamma functions has been longer than at least sixty years since the paper \cite{wendel} by J. G. Wendel was published in 1948.
\par
The motivations of bounding the ratio of two gamma functions are various, including establishment of asymptotic relation, refinements of Wallis' formula, approximation to $\pi$, needs in statistics and other mathematical sciences.
\par
In this survey paper, along one of main lines of bounding the ratio of two gamma functions, we would like to look back and analyse some known results, including Wendel's asymptotic relation, Gurland's approximation to $\pi$, Kazarinoff's refinement of Wallis' formula, Gautschi's double inequality, Watson's monotonicity, Chu's refinement of Wallis' formula, Lazarevi\'c-Lupa\c{s}'s claim on monotonic and convex properties, Kershaw's first double inequality, Elezovi\'c-Giordano-Pe\v{c}ari\'c's theorem, alternative proofs of Elezovi\'c-Giordano-Pe\v{c}ari\'c's theorem and related consequences.
\par
On the other hand, we would also like to describe some new advances in recent years on this topic, including the complete monotonicity of divided differences of the psi and polygamma functions, inequalities for sums and related results.

\section{Wendel's double inequality}

Our starting point is a paper published in 1948 by J. G. Wendel, which is the earliest related one we could search out to the best of our ability.
\par
In order to establish the classical asymptotic relation
\begin{equation}\label{wendel-approx}
\lim_{x\to\infty}\frac{\Gamma(x+s)}{x^s\Gamma(x)}=1
\end{equation}
for real $s$ and $x$, by using H\"older's inequality for integrals, J. G. Wendel~\cite{wendel} proved elegantly the double inequality
\begin{equation}\label{wendel-inequal}
\biggl(\frac{x}{x+s}\biggr)^{1-s}\le\frac{\Gamma(x+s)}{x^s\Gamma(x)}\le1
\end{equation}
for $0<s<1$ and $x>0$.

\begin{proof}[Wendel's original proof]
Let
\begin{gather*}
0<s<1,\quad p=\frac1s,\quad q=\frac{p}{p-1}=\frac1{1-s},\\
f(t)=e^{-st}t^{sx},\quad g(t)=e^{-(1-s)t}t^{(1-s)x+s-1},
\end{gather*}
and apply H\"older's inequality for integrals and the recurrent formula
\begin{equation}\label{gamma-recurrence}
\Gamma(x+1)=x\Gamma(x)
\end{equation}
for $x>0$ to obtain
\begin{equation}\label{wendel-ineq-4}
\begin{split}
\Gamma(x+s)&=\int_0^\infty e^{-t}t^{x+s-1}\td t\\
&\le\biggl(\int_0^\infty e^{-t}t^x\td t\biggr)^s\biggl(\int_0^\infty e^{-t}t^{x-1}\td t\biggr)^{1-s}\\
&=[\Gamma(x+1)]^s[\Gamma(x)]^{1-s}\\
&=x^s\Gamma(x).
\end{split}
\end{equation}
\par
Replacing $s$ by $1-s$ in~\eqref{wendel-ineq-4} we get
\begin{equation}
  \Gamma(x+1-s)\le x^{1-s}\Gamma(x),
\end{equation}
from which we obtain
\begin{equation}\label{wendel-ineq-6}
  \Gamma(x+1)\le(x+s)^{1-s}\Gamma(x+s),
\end{equation}
by substituting $x+s$ for $x$.
\par
Combining~\eqref{wendel-ineq-4} and~\eqref{wendel-ineq-6} we get
\begin{equation*}
  \frac{x}{(x+s)^{1-s}}\Gamma(x)\le\Gamma(x+s)\le x^s\Gamma(x).
\end{equation*}
Therefore, the inequality~\eqref{wendel-inequal} follows.
\par
Letting $x$ tend to infinity in~\eqref{wendel-inequal} yields~\eqref{wendel-approx} for $0<s<1$. The extension to all real $s$ is immediate on repeated application of~\eqref{gamma-recurrence}.
\end{proof}

\begin{rem}
The inequality~\eqref{wendel-inequal} can be rewritten for $0<s<1$ and $x>0$ as
\begin{equation}\label{wendel-inequal-rew}
x^{1-s}\le\frac{\Gamma(x+1)}{\Gamma(x+s)}\le(x+s)^{1-s}
\end{equation}
or
\begin{equation}\label{wendel-inequal-power}
0\le\biggl[\frac{\Gamma(x+1)}{\Gamma(x+s)}\biggr]^{1/(1-s)}-x\le s.
\end{equation}
\end{rem}

\begin{rem}
Using the recurrent formula~\eqref{gamma-recurrence} and the double inequality~\eqref{wendel-inequal-rew} repeatedly yields
\begin{equation}\label{wendel-inequal-repeat}
x^{1-s}\frac{\prod_{i=1}^m(x+i)}{\prod_{i=0}^{n-1}(x+i+s)} \le\frac{\Gamma(x+m+1)}{\Gamma(x+n+s)}\le(x+s)^{1-s} \frac{\prod_{i=1}^m(x+i)}{\prod_{i=0}^{n-1}(x+i+s)}
\end{equation}
for $x>0$ and $0<s<1$, where $m$ and $n$ are positive integers. This implies that basing on the recurrent formula~\eqref{gamma-recurrence} and the double inequality~\eqref{wendel-inequal-rew} one can bound the ratio $\frac{\Gamma(x+a)}{\Gamma(x+b)}$ for any positive numbers $x$, $a$ and $b$. Conversely, the double inequality~\eqref{wendel-inequal-repeat} reveals that one can also deduce corresponding bounds of the ratio $\frac{\Gamma(x+1)}{\Gamma(x+s)}$ for $x>0$ and $0<s<1$ from bounds of the ratio $\frac{\Gamma(x+a)}{\Gamma(x+b)}$ for positive numbers $x$, $a$ and $b$.
\end{rem}

\begin{rem}
In \cite[p.~257, 6.1.46]{abram}, the following limit was listed: For real numbers $a$ and $b$,
\begin{equation}\label{gamma-ratio-lim}
\lim_{x\to\infty}\biggl[x^{b-a}\frac{\Gamma(x+a)}{\Gamma(x+b)}\biggr]=1.
\end{equation}
The limits~\eqref{wendel-approx} and~\eqref{gamma-ratio-lim} are equivalent to each other since
\begin{equation*}
x^{t-s}\frac{\Gamma(x+s)}{\Gamma(x+t)}=\frac{\Gamma(x+s)}{x^s\Gamma(x)} \cdot\frac{x^t\Gamma(x)}{\Gamma(x+t)}.
\end{equation*}
Hence, the limit~\eqref{gamma-ratio-lim} is called as Wendel's limit in the literature of this paper.
\end{rem}

\begin{rem}
The double inequality~\eqref{wendel-inequal} or~\eqref{wendel-inequal-rew} is more meaningful than the limit~\eqref{wendel-approx} or \eqref{gamma-ratio-lim}, since the former implies the latter, but not conversely.
\end{rem}

\begin{rem}
Due to unknown reasons, Wendel's paper~\cite{wendel} and the double inequality~\eqref{wendel-inequal} or~\eqref{wendel-inequal-rew} were seemingly neglected by nearly all mathematicians for more than fifty years until it was mentioned in \cite{Merkle-JMAA-99}, to the best of my knowledge.
\end{rem}

\section{Gurland's double inequality}\label{gurland-upper-sec}

By making use of a basic theorem in mathematical statistics concerning unbiased estimators with minimum variance, J. Gurland \cite{GURLAND} presented the following inequality
\begin{equation}\label{gurland-integer-ineq}
\biggl[\frac{\Gamma((n+1)/2)}{\Gamma(n/2)}\biggr]^2<\frac{n^2}{2n+1}
\end{equation}
for $n\in\mathbb{N}$, and so taking respectively $n=2k$ and $n=2k+1$ for $k\in\mathbb{N}$ in~\eqref{gurland-integer-ineq} yields a closer approximation to $\pi$:
\begin{equation}\label{gurland-pi-ineq}
\frac{4k+3}{(2k+1)^2}\biggl[\frac{(2k)!!}{(2k-1)!!}\biggr]^2<\pi <\frac4{4k+1}\biggl[\frac{(2k)!!}{(2k-1)!!}\biggr]^2,\quad k\in\mathbb{N}.
\end{equation}

\begin{rem}
Taking respectively $n=2k$ and $n=2k-1$ for $k\in\mathbb{N}$ in~\eqref{gurland-integer-ineq} leads to
\begin{equation}\label{gurland-becomes}
\sqrt{k+\frac14}\,<\frac{\Gamma(k+1)}{\Gamma(k+1/2)}<\frac{2k}{\sqrt{4k-1}\,} =\frac{k}{\sqrt{k-1/4}\,},\quad k\in\mathbb{N}.
\end{equation}
This is better than the double inequality~\eqref{wendel-inequal-rew} for $x=k$ and $s=\frac12$.
\end{rem}

\begin{rem}
The double inequality~\eqref{gurland-pi-ineq} may be rearranged as
\begin{equation}\label{gurland-becomes-pi}
\sqrt{k+\frac14}\,<\frac{\Gamma(k+1)}{\Gamma(k+1/2)}<\frac{2k+1}{\sqrt{4k+3}\,} =\frac{k+1/2}{\sqrt{k+1/2+1/4}\,},\quad k\in\mathbb{N}.
\end{equation}
\par
It is easy to see that the upper bound in~\eqref{gurland-becomes-pi} is better than the corresponding one in~\eqref{gurland-becomes}. This phenomenon seemingly hints us that sharper bounds for the ratio $\frac{\Gamma(k+1)}{\Gamma(k+1/2)}$ can be obtained only if letting $m\in\mathbb{N}$ in $n=2m-1$ is larger in~\eqref{gurland-integer-ineq}. However, this is an illusion, since the lower bound of the following double inequality
\begin{multline}\label{frac12-gen-ineq}
\frac{k^2}{2k+1}\cdot\frac{\sqrt{4(k+m)-3}\,}{k+m-1}\prod_{i=0}^{m-1}\biggl[1+\frac1{2(k+i)}\biggr] <\frac{\Gamma(k+1)}{\Gamma(k+1/2)}\\
<\biggl(\frac{2k}{2k-1}\biggr)^2\frac{2(k+m)-3}{\sqrt{4(k+m)-5}\,} \prod_{i=0}^{m-1}\biggl[1-\frac1{2(k+i)}\biggr],
\end{multline}
which is derived from taking respectively $n=2(k+m-1)$ and $n=2(k+m-1)-1$ for $k\in\mathbb{N}$ in~\eqref{gurland-integer-ineq}, is decreasing and the upper bound of it is increasing with respect to $m$. Then how to explain the occurrence that the upper bound in~\eqref{gurland-becomes-pi} is stronger than the corresponding one in~\eqref{gurland-becomes}?
\end{rem}

\begin{rem}
The left-hand side inequality in~\eqref{gurland-becomes} or~\eqref{gurland-becomes-pi} may be rearranged as
\begin{equation}\label{left-ineq-rew}
\frac14<\biggl[\frac{\Gamma(k+1)}{\Gamma(k+1/2)}\biggr]^2-k<\frac14+\frac1{4(4k+3)},\quad k\in\mathbb{N}.
\end{equation}
From this, it is easier to see that the inequality~\eqref{gurland-integer-ineq} refines the double inequality~\eqref{wendel-inequal-rew} for $x=k$ and $s=\frac12$.
\end{rem}

\begin{rem}
It is noted that the inequality~\eqref{gurland-integer-ineq} was recovered in \cite{Chu-monthly} and extended in \cite{Kazarinoff-56} by different approaches respectively. See Section~\ref{Wallis-section} and Section~\ref{chu-sec} below.
\end{rem}

\begin{rem}
Just like the paper \cite{wendel}, Gurland's paper \cite{GURLAND} was ignored except it was mentioned in \cite{dutka, Rao-Uppuluri}. The famous monograph \cite{mit} recorded neither of the papers \cite{GURLAND, wendel}.
\end{rem}

\section{Kazarinoff's double inequality}\label{Wallis-section}

Starting from one form of the celebrated formula of John Wallis:
\begin{equation}\label{John-Wallis-ineq}
\frac1{\sqrt{\pi(n+1/2)}}<\frac{(2n-1)!!}{(2n)!!}<\frac1{\sqrt{\pi n}},\quad n\in\mathbb{N},
\end{equation}
which had been quoted for more than a century before 1950s by writers of textbooks, D. K. Kazarinoff proved in~\cite{Kazarinoff-56} that the sequence $\theta(n)$ defined by
\begin{equation}\label{theta-dfn-kazar}
\frac{(2n-1)!!}{(2n)!!} =\frac1{\sqrt{\pi[n+\theta(n)]}}
\end{equation}
satisfies $\frac14<\theta(n)<\frac12$ for $n\in\mathbb{N}$. This implies
\begin{equation}\label{Wallis'inequality}
\frac1{\sqrt{\pi(n+1/2)}}<\frac{(2n-1)!!}{(2n)!!}
<\frac1{\sqrt{\pi(n+1/4)}},\quad n\in\mathbb{N}.
\end{equation}

\begin{rem}
It was said in \cite{Kazarinoff-56} that it is unquestionable that inequalities similar to~\eqref{Wallis'inequality} can be improved indefinitely but at a sacrifice of simplicity, which is why the inequality~\eqref{John-Wallis-ineq} had survived so long.
\end{rem}

\begin{rem}
Kazarinoff's proof of~\eqref{Wallis'inequality} is based upon the property
\begin{equation}\label{Phi-ineq}
[\ln\phi(t)]''-\{[\ln\phi(t)]'\}^2>0
\end{equation}
of the function
\begin{equation}
\phi(t)=\int_0^{\pi/2}\sin^tx\td x=\frac{\sqrt\pi\,}2\cdot\frac{\Gamma((t+1)/2)}{\Gamma((t+2)/2)}
\end{equation}
for $-1<t<\infty$. The inequality~\eqref{Phi-ineq} was proved by making use of the well-known Legendre's formula
\begin{equation}\label{Legendre's-formula}
\psi(x)=-\gamma+\int_0^1\frac{t^{x-1}-1}{t-1}\td t
\end{equation}
for $x>0$ and estimating the integrals
\begin{equation}
\int_0^1\frac{x^t}{1+x}\td x\quad\text{and}\quad \int_0^1\frac{x^t\ln x}{1+x}\td x.
\end{equation}
Since~\eqref{Phi-ineq} is equivalent to the statement that the reciprocal of $\phi(t)$ has an everywhere negative second derivative, therefore, for any positive $t$, $\phi(t)$ is less than the harmonic mean of $\phi(t-1)$ and $\phi(t+1)$, which implies
\begin{equation}\label{karz-2.17-ineq}
\frac{\Gamma((t+1)/2)}{\Gamma((t+2)/2)}<\frac2{\sqrt{2t+1}},\quad t>-\frac12.
\end{equation}
As a subcase of this result, the right-hand side inequality in~\eqref{Wallis'inequality} is established.
\end{rem}

\begin{rem}
Using the recurrent formula~\eqref{gamma-recurrence} in~\eqref{karz-2.17-ineq} gives
\begin{equation}\label{gurland-integer-ineq-ext}
\biggl[\frac{\Gamma((t+1)/2)}{\Gamma(t/2)}\biggr]^2<\frac{t^2}{2t+1}
\end{equation}
for $t>0$, which extends the inequality~\eqref{gurland-integer-ineq}. This shows that Kazarinoff's paper \cite{Kazarinoff-56} contains much general conclusions and that all results in \cite{GURLAND} stated in Section~\ref{gurland-upper-sec} are consequences of the inequality~\eqref{gurland-integer-ineq-ext}, as showed below.
\par
Replacing $t$ by $2t$ in~\eqref{karz-2.17-ineq} or~\eqref{gurland-integer-ineq-ext} and rearranging yield
\begin{equation}
\frac{\Gamma(t+1)}{\Gamma(t+1/2)}>\sqrt{t+\frac14} \quad \Longleftrightarrow\quad \biggl[\frac{\Gamma(t+1)}{\Gamma(t+1/2)}\biggr]^2-t>\frac14
\end{equation}
for $t>0$, which extends the left-hand side inequality in~\eqref{gurland-becomes} and~\eqref{gurland-becomes-pi}. Replacing $t$ by $2t-1$ in~\eqref{karz-2.17-ineq} or~\eqref{gurland-integer-ineq-ext} produces
\begin{equation}
\frac{\Gamma(t+1)}{\Gamma(t+1/2)}<\frac{2t}{\sqrt{4t-1}\,}
\end{equation}
for $t>\frac12$, which extends the right-hand side inequality in~\eqref{gurland-becomes}. Replacing $t$ by $2t+1$ in~\eqref{karz-2.17-ineq} or~\eqref{gurland-integer-ineq-ext} and rearranging gives
\begin{equation}
\frac{\Gamma(t+1)}{\Gamma(t+1/2)}<\frac{2t+1}{\sqrt{4t+3}\,}
\end{equation}
for $t>-\frac12$, which extends the right-hand side inequality in~\eqref{gurland-becomes-pi}.
\end{rem}

\begin{rem}
By the well-known Wallis cosine formula \cite{WallisFormula.html}, the sequence $\theta(n)$ defined by \eqref{theta-dfn-kazar} may be rearranged as
\begin{equation}\label{theta(n)}
\theta(n)=\frac1\pi\biggl[\frac{(2n)!!}{(2n-1)!!}\biggr]^2-n =\biggl[\frac{\Gamma(n+1)}{\Gamma(n+1/2)}\biggr]^2-n
\end{equation}
for $n\in\mathbb{N}$. Then the inequality~\eqref{Wallis'inequality} is equivalent to
\begin{equation}
\frac14<\biggl[\frac{\Gamma(n+1)}{\Gamma(n+1/2)}\biggr]^2-n<\frac12,\quad n\in\mathbb{N}.
\end{equation}
\end{rem}

\begin{rem}\label{kazarinoff-rem-5}
The inequality~\eqref{Phi-ineq} may be rewritten as
\begin{equation}
\psi'\biggl(\frac{t+1}2\biggr)-\psi'\biggl(\frac{t+2}2\biggr) >\biggl[\psi\biggl(\frac{t+1}2\biggr)-\psi\biggl(\frac{t+2}2\biggr)\biggr]^2
\end{equation}
for $t>-1$. Letting $u=\frac{t+1}2$ in the above inequality yields
\begin{equation}
\psi'(u)-\psi'\biggl(u+\frac12\biggr) >\biggl[\psi(u)-\psi\biggl(u+\frac12\biggr)\biggr]^2
\end{equation}
for $u>0$. This inequality has been generalized in \cite{Comp-Mon-Digamma-Trigamma-Divided.tex} to the complete monotonicity of a function involving divided differences of the digamma and trigamma functions as follows.

\begin{thm}\label{CMDT-divided-thm}
For real numbers $s$, $t$, $\alpha=\min\{s,t\}$ and $\lambda$, let
\begin{equation}\label{Delta-lambda-dfn}
\Delta_{s,t;\lambda}(x)=\begin{cases}\bigg[\dfrac{\psi(x+t) -\psi(x+s)}{t-s}\bigg]^2
+\lambda\dfrac{\psi'(x+t)-\psi'(x+s)}{t-s},&s\ne t\\
[\psi'(x+s)]^2+\lambda\psi''(x+s),&s=t
\end{cases}
\end{equation}
on $(-\alpha,\infty)$. Then the function $\Delta_{s,t;\lambda}(x)$ has the following complete monotonicity:
\begin{enumerate}
\item
For $0<|t-s|<1$,
\begin{enumerate}
\item
the function $\Delta_{s,t;\lambda}(x)$ is completely monotonic on $(-\alpha,\infty)$ if and only if $\lambda\le1$,
\item
so is the function $-\Delta_{s,t;\lambda}(x)$ if and only if $\lambda\ge\frac1{|t-s|}$;
\end{enumerate}
\item
For $|t-s|>1$,
\begin{enumerate}
\item
the function $\Delta_{s,t;\lambda}(x)$ is completely monotonic on $(-\alpha,\infty)$ if and only if $\lambda\le\frac1{|t-s|}$,
\item
so is the function $-\Delta_{s,t;\lambda}(x)$ if and only if $\lambda\ge1$;
\end{enumerate}
\item
For $s=t$, the function $\Delta_{s,s;\lambda}(x)$ is completely monotonic on $(-s,\infty)$ if and only if $\lambda\le1$;
\item
For $|t-s|=1$,
\begin{enumerate}
\item
the function $\Delta_{s,t;\lambda}(x)$ is completely monotonic if and only if $\lambda<1$,
\item
so is the function $-\Delta_{s,t;\lambda}(x)$ if and only if $\lambda>1$,
\item
and $\Delta_{s,t;1}(x)\equiv0$.
\end{enumerate}
\end{enumerate}
\end{thm}
Taking in Theorem~\ref{CMDT-divided-thm} $\lambda=s-t>0$ produces that the function $\frac{\Gamma(x+s)}{\Gamma(x+t)}$ on $(-t,\infty)$ is increasingly convex if $s-t>1$ and increasingly concave if $0<s-t<1$.
\end{rem}

\section{Watson's monotonicity}\label{Watson-sec}

In 1959, motivated by the result in \cite{Kazarinoff-56} mentioned in Section~\ref{Wallis-section}, G.~N.~Watson~\cite{waston} observed that
\begin{multline}\label{watson-formula}
\frac1x\cdot\frac{[\Gamma(x+1)]^2}{[\Gamma(x+1/2)]^2} ={}_2F_1\biggl(-\frac12,-\frac12;x;1\biggr)\\*
=1+\frac1{4x}+\frac1{32x(x+1)} +\sum_{r=3}^\infty\frac{[(-1/2)\cdot(1/2)\cdot(3/2)\cdot(r-3/2)]^2} {r!x(x+1)\dotsm(x+r-1)}
\end{multline}
for $x>-\frac12$, which implies the much general function
\begin{equation}\label{theta-dfn}
\theta(x)=\biggl[\frac{\Gamma(x+1)}{\Gamma(x+1/2)}\biggr]^2-x
\end{equation}
for $x>-\frac12$, whose special case is the sequence $\theta(n)$ for $n\in\mathbb{N}$ defined in~\eqref{theta-dfn-kazar} or~\eqref{theta(n)}, is decreasing and
\begin{equation}
\lim_{x\to\infty}\theta(x)=\frac14\quad \text{and}\quad \lim_{x\to(-1/2)^+}\theta(x)=\frac12.
\end{equation}
This apparently implies the sharp inequalities
\begin{equation}\label{theta-l-u-b}
\frac14<\theta(x)<\frac12
\end{equation}
for $x>-\frac12$,
\begin{equation}\label{watson-special-ineq}
\sqrt{x+\frac14}\,< \frac{\Gamma(x+1)}{\Gamma(x+1/2)}\le \sqrt{x+\frac14+\biggl[\frac{\Gamma(3/4)}{\Gamma(1/4)}\biggr]^2}\, =\sqrt{x+0.36423\dotsm}
\end{equation}
for $x\ge-\frac14$, and, by Wallis cosine formula~\cite{WallisFormula.html},
\begin{equation}\label{best-bounds-Wallis}
\frac{1}{\sqrt{\pi(n+{4}/{\pi}-1)}}\le\frac{(2n-1)!!}{(2n)!!}
<\frac{1}{\sqrt{\pi(n+1/4)}},\quad n\in\mathbb{N}.
\end{equation}
\par
In \cite{waston}, an alternative proof of the double inequality~\eqref{theta-l-u-b} was also provided as follows: Let
\begin{equation}
f(x)=\frac2{\sqrt\pi\,}\int_0^{\pi/2}\cos^{2x}t\td t =\frac2{\sqrt\pi\,}\int_0^\infty\exp\bigl(-xt^2\bigr) \frac{t\exp(-t^2/2)}{\sqrt{1-\exp(-t^2)}\,}\td t
\end{equation}
for $x>\frac12$. By using the fairly obvious inequalities
\begin{equation}
\sqrt{1-\exp\bigl(-t^2\bigr)}\,\le t
\end{equation}
and
\begin{equation}
\frac{t\exp(-t^2/4)}{\sqrt{1-\exp(-t^2)}\,}=\frac{t}{\sqrt{2\sinh(t^2/2)}}\le1,
\end{equation}
we have, for $x>-\frac14$,
\begin{equation*}
\frac1{\sqrt{\pi}\,}\int_0^\infty\exp\bigl(-(x+1/2)t^2\bigr)\td t<f(x) <\frac1{\sqrt{\pi}\,}\int_0^\infty\exp\bigl(-(x+1/4)t^2\bigr)\td t,
\end{equation*}
that is to say
\begin{equation}
\frac1{\sqrt{x+1/2}\,}<f(x)<\frac1{\sqrt{x+1/4}\,}.
\end{equation}

\begin{rem}
It is easy to see that the inequality~\eqref{watson-special-ineq} extends and improves~\eqref{wendel-inequal-rew} if $s=\frac12$, say nothing of~\eqref{gaut-ineq-1} and~\eqref{gaut-ineq-2} if $s=\frac12$.
\end{rem}

\begin{rem}
The left-hand side inequality in~\eqref{best-bounds-Wallis} is better than the corresponding one in~\eqref{Wallis'inequality} but worse than the corresponding one in~\eqref{gurland-pi-ineq} for $n\ge2$.
\end{rem}

\begin{rem}
The formula~\eqref{watson-formula} implies the complete monotonicity of the function $\theta(x)$ defined by~\eqref{theta-dfn} on $\bigl(-\frac12,\infty\bigr)$.
\end{rem}

\section{Gautschi's double inequalities}

The main aim of the paper \cite{gaut} was to establish the double inequality
\begin{equation}\label{gaut-3-ineq}
\frac{(x^p+2)^{1/p}-x}2<e^{x^p}\int_x^\infty e^{-t^p}\td t\le c_p\biggl[\biggl(x^p+\frac1{c_p}\biggr)^{1/p}-x\biggr]
\end{equation}
for $x\ge0$ and $p>1$, where
\begin{equation}
c_p=\biggl[\Gamma\biggl(1+\frac1p\biggr)\biggr]^{p/(p-1)}
\end{equation}
or $c_p=1$.
\par
By an easy transformation, the inequality~\eqref{gaut-3-ineq} was written in terms of the complementary gamma function
\begin{equation}
\Gamma(a,x)=\int_x^\infty e^{-t}t^{a-1}\td t
\end{equation}
as
\begin{equation}\label{gaut-4-ineq}
\frac{p[(x+2)^{1/p}-x^{1/p}]}2<e^x\Gamma\biggl(\frac1p,x\biggr)\le pc_p\biggl[\biggl(x+\frac1{c_p}\biggr)^{1/p}-x^{1/p}\biggr]
\end{equation}
for $x\ge0$ and $p>1$. In particular, if letting $p\to\infty$, the double inequality
\begin{equation}
\frac12\ln\biggl(1+\frac2x\biggr)\le e^xE_1(x)\le\ln\biggl(1+\frac1x\biggr)
\end{equation}
for the exponential integral $E_1(x)=\Gamma(0,x)$ for $x>0$ was derived from~\eqref{gaut-4-ineq}, in which the bounds exhibit the logarithmic singularity of $E_1(x)$ at $x=0$.
\par
As a direct consequence of the inequality~\eqref{gaut-4-ineq} for $p=\frac1s$, $x=0$ and $c_p=1$, the following simple inequality for the gamma function was deduced:
\begin{equation}\label{gaut-none-ineq}
2^{s-1}\le\Gamma(1+s)\le1,\quad 0\le s\le 1.
\end{equation}
The second main result of the paper \cite{gaut} was a sharper and more general inequality
\begin{equation}\label{gaut-6-ineq}
e^{(s-1)\psi(n+1)}\le\frac{\Gamma(n+s)}{\Gamma(n+1)}\le n^{s-1}
\end{equation}
for $0\le s\le1$ and $n\in\mathbb{N}$ than~\eqref{gaut-none-ineq} by proving that the function
\begin{equation}
f(s)=\frac1{1-s}\ln\frac{\Gamma(n+s)}{\Gamma(n+1)}
\end{equation}
is monotonically decreasing for $0\le s<1$. Since $\psi(n)<\ln n$, it was derived from the inequality~\eqref{gaut-6-ineq} that
\begin{equation}\label{gaut-6-ineq-simp}
\biggl(\frac1{n+1}\biggr)^{1-s}\le\frac{\Gamma(n+s)}{\Gamma(n+1)}\le\biggl(\frac1n\biggr)^{1-s}, \quad 0\le s\le1,
\end{equation}
which was also rewritten as
\begin{equation}\label{euler-gaut}
\frac{n!(n+1)^{s-1}}{(s+1)(s+2)\dotsm(s+n-1)}\le\Gamma(1+s) \le\frac{(n-1)!n^s}{(s+1)(s+2)\dotsm(s+n-1)}
\end{equation}
and so a simple proof of Euler's product formula in the segment $0\le s\le1$ was showed by letting $n\to\infty$ in~\eqref{euler-gaut}.

\begin{rem}
The double inequalities~\eqref{gaut-6-ineq} and~\eqref{gaut-6-ineq-simp} can be further rearranged as
\begin{equation}\label{gaut-ineq-1}
n^{1-s}\le\frac{\Gamma(n+1)}{\Gamma(n+s)}\le\exp((1-s)\psi(n+1))
\end{equation}
and
\begin{equation}\label{gaut-ineq-2}
n^{1-s}\le\frac{\Gamma(n+1)}{\Gamma(n+s)}\le (n+1)^{1-s}
\end{equation}
for $n\in\mathbb{N}$ and $0\le s\le 1$.
\end{rem}

\begin{rem}
The upper bounds in~\eqref{wendel-inequal-rew} and~\eqref{gaut-ineq-1} have the following relationship
\begin{equation}\label{wendel-gautschi-comp}
(n+s)^{1-s}\le\exp((1-s)\psi(n+1))
\end{equation}
for $0\le s\le\frac12$ and $n\in\mathbb{N}$, and the inequality~\eqref{wendel-gautschi-comp} reverses for $s>e^{1-\gamma}-1=0.52620\dotsm$, since the function
\begin{equation}\label{Q(x)-dfn}
Q(x)=e^{\psi(x+1)}-x
\end{equation}
was proved in~\cite[Theorem~2]{Infinite-family-Digamma.tex} to be strictly decreasing on $(-1,\infty)$, with
\begin{equation}\label{Q-infty-lim}
\lim_{x\to\infty}Q(x)=\frac12.
\end{equation}
This means that Wendel's double inequality~\eqref{wendel-inequal-rew} and Gautschi's first double inequality~\eqref{gaut-ineq-1} are not included each other but they all contain Gautschi's second double inequality~\eqref{gaut-ineq-2}.
\end{rem}

\begin{rem}
The right-hand side inequality in \eqref{gaut-ineq-1} may be rearranged as
\begin{equation}
\biggl[\frac{\Gamma(n+1)}{\Gamma(n+s)}\biggr]^{1/(1-s)}\le\exp(\psi(n+1)),\quad n\in\mathbb{N}.
\end{equation}
This suggests us the following double inequality
\begin{equation}\label{aplha-beta-ineq}
\exp(\psi(\alpha(x)))<\biggl[\frac{\Gamma(x+t)}{\Gamma(x+s)}\biggr]^{1/(t-s)}\le\exp(\psi(\beta(x)))
\end{equation}
for real numbers $s$, $t$ and $x\in(-\min\{s,t\},\infty)$, where $\alpha(x)\sim x$ and $\beta(x)\sim x$ as $x\to\infty$. For detailed information on the type of inequalities like \eqref{aplha-beta-ineq}, please refer to \cite{bounds-two-gammas.tex} and related references therein.
\end{rem}

\begin{rem}
The inequality~\eqref{gaut-ineq-2} can be rewritten as
\begin{equation}\label{gaut-ineq-2-rew}
0\le\biggl[\frac{\Gamma(n+1)}{\Gamma(n+s)}\biggr]^{1/(1-s)}-n\le1
\end{equation}
for $n\in\mathbb{N}$ and $0\le s\le 1$.
\end{rem}

\begin{rem}
In the texts of the reviews on the paper~\cite{gaut} by the Mathematical Reviews and the Zentralblatt MATH, there is no a word to comment on inequalities in~\eqref{gaut-ineq-1} and~\eqref{gaut-ineq-2}. However, these two double inequalities later became a major source of a series of study on bounding the ratio of two gamma functions.
\end{rem}

\section{Chu's double inequality}\label{chu-sec}

In 1962, by discussing that
\begin{equation}\label{chu-proof-idea}
b_{n+1}(c)\gtreqqless b_n(c)
\end{equation}
if and only if $(1-4c)n+1-3c\gtreqqless0$, where
\begin{equation}\label{b-n(c)}
b_n(c)=\frac{(2n-1)!!}{(2n)!!}\sqrt{n+c}\,,
\end{equation}
it was demonstrated in~\cite[Theorem~1]{Chu-monthly} that
\begin{equation}\label{Chu-bounds-wallis}
\frac1{\sqrt{\pi[n+(n+1)/(4n+3)]}\,}<\frac{(2n-1)!!}{(2n)!!}<\frac1{\sqrt{\pi(n+1/4)}},\quad n\in\mathbb{N}.
\end{equation}
As an application of~\eqref{Chu-bounds-wallis}, by using $\Gamma\bigl(\frac12\bigr)=\sqrt\pi$  and~\eqref{gamma-recurrence}, the following double inequality
\begin{equation}\label{chu-ineq-ratio}
\sqrt{\frac{2n-3}4}\,<\frac{\Gamma(n/2)}{\Gamma(n/2-1/2)} \le\sqrt{\frac{(n-1)^2}{2n-1}}\,
\end{equation}
for positive integers $n\ge2$ was given in~\cite[Theorem~2]{Chu-monthly}.

\begin{rem}
After letting $n=2k+1$ the inequality~\eqref{chu-ineq-ratio} becomes
\begin{equation}\label{chu-becomes}
\sqrt{k-\frac14}\,<\frac{\Gamma(k+1/2)}{\Gamma(k)}<\frac{k}{\sqrt{k+1/4}},
\end{equation}
which is same as \eqref{gurland-becomes}. Taking $n=2k+2$ in leads to inequalities~\eqref{gurland-becomes-pi} and \eqref{left-ineq-rew}.
\par
Notice that the reasoning directions in the two papers \cite{Chu-monthly, GURLAND} are opposite:
\begin{equation}
\frac{(2n-1)!!}{(2n)!!}\quad \begin{gathered} \overset{\text{\cite{Chu-monthly}}}{\Longrightarrow} \\[-0.6em] \underset{\text{\cite{GURLAND}}}{\Longleftarrow}\end{gathered} \quad \frac{\Gamma(n/2)}{\Gamma(n/2-1/2)}.
\end{equation}
To some extent, the results obtained by Gurland in~\cite{GURLAND} and by Chu in~\cite{Chu-monthly} are equivalent to each other and they are all special cases of those obtained by Kazarinoff in~\cite{Kazarinoff-56}.
\end{rem}

\begin{rem}
By Wallis cosine formula~\cite{WallisFormula.html}, the sequence \eqref{b-n(c)} may be rewritten as
\begin{equation}\label{b-n(c)-rew}
b_n(c)=\frac1{\sqrt\pi\,}\cdot\frac{\Gamma(n+1/2)}{\Gamma(n+1)}\sqrt{n+c}\, \triangleq\frac1{\sqrt\pi\,}B_c(n)
\end{equation}
for $n\in\mathbb{N}$. Therefore, Chu discussed equivalently the necessary and sufficient conditions such that the sequence $B_c(n)$ for $n\in\mathbb{N}$ is monotonic.
\par
Recently, necessary and sufficient conditions for the general function
\begin{equation}\label{h-def-sandor-new}
H_{a,b,c}(x)=(x+c)^{b-a}\frac{\Gamma(x+a)}{\Gamma(x+b)}
\end{equation}
on $(-\rho,\infty)$, where $a$, $b$ and $c$ are real numbers and $\rho=\min\{a,b,c\}$, to be logarithmically completely monotonic are presented in \cite{sandor-gamma-3-note.tex-final, sandor-gamma-3-note.tex}. A positive function $f$ is said to be logarithmically completely monotonic on an interval $I\subseteq\mathbb{R}$ if it has derivatives of all orders on $I$ and its logarithm $\ln f$ satisfies $(-1)^k[\ln f(x)]^{(k)}\ge0$ for $k\in\mathbb{N}$ on $I$, see~\cite{Atanassov, CBerg, minus-one}.
\end{rem}

\section{Lazarevi\'c-Lupa\c{s}'s claim}\label{lazarevic-sec}
In 1974, among other things, the function
\begin{equation}\label{Lazarevic-Lupas-funct}
\theta_{\alpha}(x)=\biggl[\frac{\Gamma(x+1)}{\Gamma(x+\alpha)}\biggr]^{1/(1-\alpha)}-x
\end{equation}
on $(0,\infty)$ for $\alpha\in(0,1)$ was claimed in~\cite[Theorem~2]{Lazarevic} to be decreasing and convex, and so
\begin{equation}
\frac{\alpha}2<\biggl[\frac{\Gamma(x+1)}{\Gamma(x+\alpha)}\biggr]^{1/(1-\alpha)}-x \le[\Gamma(\alpha)]^{1/(1-\alpha)}.
\end{equation}

\begin{rem}
The proof of \cite[Theorem~2]{Lazarevic} is wrong, see \cite[Remark~3.3]{Alzer-Math-Nachr-2001} and \cite[p.~240]{egp}. However, the statements in \cite[Theorem~2]{Lazarevic} are correct and this is the first time to try to investigate the monotonic and convex properties of the much general function $\theta_{\alpha}(x)$.
\end{rem}

\section{Kershaw's first double inequality}\label{kershaw-sec}

In 1983, motivated by the inequality~\eqref{gaut-ineq-2} obtained in~\cite{gaut}, among other things, Kershaw presented in~\cite{kershaw} the following double inequality
\begin{equation}\label{gki1}
\biggl(x+\frac{s}2\biggr)^{1-s}<\frac{\Gamma(x+1)}{\Gamma(x+s)}
<\biggl[x-\frac12+\biggl(s+\frac14\biggr)^{1/2}\biggr]^{1-s}
\end{equation}
for $0<s<1$ and $x>0$. In the literature, it is called as Kershaw's first double inequality for the ratio of two gamma functions.

\begin{proof}[Kershaw's proof for~\eqref{gki1}]
Define the function $g_\beta$ by
\begin{equation}\label{kershaw-g-dfn}
  g_\beta(x)=\frac{\Gamma(x+1)}{\Gamma(x+s)}(x+\beta)^{s-1}
\end{equation}
for $x>0$ and $0<s<1$, where the parameter $\beta$ is to be determined.
\par
It is not difficult to show, with the aid of Wendel's limit~\eqref{gamma-ratio-lim}, that
\begin{equation}\label{kershaw-2.3}
\lim_{x\to\infty}g_\beta(x)=1.
\end{equation}
\par
To prove the double inequality~\eqref{gki1} define
\begin{equation}\label{G(x)-dfn}
  G(x)=\frac{g_\beta(x)}{g_\beta(x+1)}=\frac{x+s}{x+1}\biggl(\frac{x+\beta+1}{x+\beta}\biggr)^{1-s},
\end{equation}
from which it follows that
\begin{equation*}
  \frac{G'(x)}{G(x)}=\frac{(1-s)[(\beta^2+\beta-s)+(2\beta-s)x]}{(x+1)(x+s)(x+\beta)(x+\beta+1)}.
\end{equation*}
This will leads to
\begin{enumerate}
  \item
  if $\beta=\frac{s}2$, then $G'(x)<0$ for $x>0$;
  \item
  if $\beta=-\frac12+\bigl(s+\frac14\bigr)^{1/2}$, then $G'(x)>0$ for $x>0$.
\end{enumerate}
\par
Consequently if $\beta=\frac{s}2$ then $G$ strictly decreases, and since $G(x)\to1$ as $x\to\infty$ it follows that $G(x)>1$ for $x>0$. But, from~\eqref{kershaw-2.3}, this implies that $g_\beta(x)>g_\beta(x+1)$ for $x>0$, and so $g_\beta(x)>g_\beta(x+n)$. Take the limit as $n\to\infty$ to give the result that $g_\beta(x)>1$, which can be rewritten as the left-hand side inequality in~\eqref{gki1}. The corresponding upper bound can be verified by a similar argument when $\beta=-\frac12+\bigl(s+\frac14\bigr)^{1/2}$, the only difference being that in this case $g_\beta$ strictly increases to unity.
\end{proof}

\begin{rem}
The spirit of Kershaw's proof is similar to Chu's in \cite[Theorem~1]{Chu-monthly}, as showed by \eqref{chu-proof-idea}. This idea or method was also utilized independently in \cite{Giordano-Laforgia-01-jcam, GLP-98, unify-kershaw, laforgia-mc-1984, lorch-ultra, Palumbo-98-jcam} to construct for various purposes a number of inequalities of the type
\begin{equation}
(x+\alpha)^{s-1}<\frac{\Gamma(x+s)}{\Gamma(x+1)}<(x+\beta)^{s-1}
\end{equation}
for $s>0$ and real number $x\ge0$.
\end{rem}

\begin{rem}
It is easy to see that the inequality~\eqref{gki1} refines and extends the inequality~\eqref{wendel-inequal-rew}, say nothing of~\eqref{gaut-ineq-2}.
\end{rem}

\begin{rem}
The inequality~\eqref{gki1} may be rearranged as
\begin{equation}\label{gki1-rew}
\frac{s}2<\biggl[\frac{\Gamma(x+1)}{\Gamma(x+s)}\biggr]^{1/(1-s)}-x <\biggl(s+\frac14\biggr)^{1/2}-\frac12
\end{equation}
for $x>0$ and $0<s<1$.
\end{rem}

\section{Elezovi\'c-Giordano-Pe\v{c}ari\'c's theorem}

The inequalities~\eqref{wendel-inequal-power}, \eqref{left-ineq-rew}, \eqref{gaut-ineq-2-rew} and~\eqref{gki1-rew}, the sequence~\eqref{theta(n)} and the function~\eqref{theta-dfn} and~\eqref{Lazarevic-Lupas-funct} strongly suggest us to consider the monotonic and convex properties of the general function
\begin{equation}
z_{s,t}(x)=\begin{cases}
\bigg[\dfrac{\Gamma(x+t)}{\Gamma(x+s)}\bigg]^{1/(t-s)}-x,&s\ne t\\
e^{\psi(x+s)}-x,&s=t
\end{cases}
\end{equation}
for $x\in(-\alpha,\infty)$, where $s$ and $t$ are two real numbers and $\alpha=\min\{s,t\}$.
\par
In 2000, N. Elezovi\'c, C. Giordano and J. Pe\v{c}ari\'c gave in \cite[Theorem~1]{egp} a perfect solution to the monotonic and convex properties of the function $z_{s,t}(x)$ as follows.

\begin{thm}\label{egp-thm-mon}
The function $z_{s,t}(x)$ is either convex and decreasing for $|t-s|<1$ or concave and increasing for $|t-s|>1$.
\end{thm}

\begin{rem}
Direct computation yields
\begin{equation}\label{z-2-deriv}
\frac{z''_{s,t}(x)}{z_{s,t}(x)+x}=\biggl[\frac{\psi(x+t)-\psi(x+s)}{t-s}\biggr]^2 +\frac{\psi'(x+t)-\psi'(x+s)}{t-s}.
\end{equation}
To prove the positivity of the function~\eqref{z-2-deriv}, the following formula and inequality are used as basic tools in the proof of \cite[Theorem~1]{egp}:
\begin{enumerate}
\item
For $x>-1$,
\begin{equation}
\psi(x+1)=-\gamma+\sum_{k=1}^\infty\biggl(\frac1k-\frac1{x+k}\biggr).
\end{equation}
\item
If $a\le b<c\le d$, then
\begin{equation}
\frac1{ab}+\frac1{cd}>\frac1{ac}+\frac1{bd}.
\end{equation}
\end{enumerate}
\end{rem}

\begin{rem}
As consequences of Theorem~\ref{egp-thm-mon}, the following useful conclusions are derived:
\begin{enumerate}
\item
The function
\begin{equation}\label{e-psi-x}
e^{\psi(x+t)}-x
\end{equation}
for all $t>0$ is decreasing and convex from $(0,\infty)$ onto $\bigl(e^{\psi(t)},t-\frac12\bigr)$.
\item
For all $x>0$,
\begin{equation}\label{egp-ineq-psi-der}
\psi'(x)e^{\psi(x)}<1.
\end{equation}
\item
For all $x>0$ and $t>0$,
\begin{equation}\label{ineq-3.52}
\ln\biggl(x+t-\frac12\biggr)<\psi(x+t)<\ln\Bigl(x+e^{\psi(t)}\Bigr).
\end{equation}
\item
For $x>-\alpha$, the inequality
\begin{equation}
\bigg[\dfrac{\Gamma(x+t)}{\Gamma(x+s)}\bigg]^{1/(t-s)}<\frac{t-s}{\psi(x+t)-\psi(x+s)}
\end{equation}
holds if $|t-s|<1$ and reverses if $|t-s|>1$.
\end{enumerate}
\end{rem}

\begin{rem}\label{rem-3.19.6}
In fact, the function~\eqref{e-psi-x} is deceasing and convex on $(-t,\infty)$ for all $t\in\mathbb{R}$. See \cite[Theorem~2]{Infinite-family-Digamma.tex}.
\end{rem}

\begin{rem}\label{rem-3.19.7}
It is clear that the double inequality~\eqref{ineq-3.52} can be deduced directly from the decreasingly monotonic property of~\eqref{e-psi-x}. Furthermore, from the decreasingly monotonic and convex properties of~\eqref{e-psi-x} on $(-t,\infty)$, the inequality~\eqref{egp-ineq-psi-der} and
\begin{equation}\label{positivity}
\psi''(x)+[\psi'(x)]^2>0
\end{equation}
on $(0,\infty)$ can be derived straightforwardly.
\end{rem}

\section{Recent advances}

Finally, we would like to state some new results related to or originated from Elezovi\'c-Giordano-Pe\v{c}ari\'c's Theorem~\ref{egp-thm-mon} above.

\subsection{Alternative proofs of Elezovi\'c-Giordano-Pe\v{c}ari\'c's theorem}\label{proofs-qi}

The key step of verifying Theorem~\ref{egp-thm-mon} is to prove the positivity of the right-hand side in~\eqref{z-2-deriv} in which involves divided differences of the digamma and trigamma functions. The biggest barrier or difficulty to prove the positivity of~\eqref{z-2-deriv} is mainly how to deal with the squared term in~\eqref{z-2-deriv}.

\subsubsection{Chen's proof}
In \cite{201-05-JIPAM}, the barrier mentioned above was overcome by virtute of the well-known convolution theorem \cite{LaplaceTransform.html} for Laplace transforms and so Theorem~\ref{egp-thm-mon} for the special case $s+1>t>s\ge0$ was proved. Perhaps this is the first try to provide an alternative of Theorem~\ref{egp-thm-mon}, although it was partially successful formally.

\subsubsection{Qi-Guo-Chen's proof}
For real numbers $\alpha$ and $\beta$ with $(\alpha,\beta)\not\in\{(0,1),(1,0)\}$ and $\alpha\ne\beta$, let
\begin{equation}\label{q-dfn}
q_{\alpha,\beta}(t)=
\begin{cases}
\dfrac{e^{-\alpha t}-e^{-\beta t}}{1-e^{-t}},&t\ne0,\\
\beta-\alpha,&t=0.
\end{cases}
\end{equation}
In \cite{notes-best.tex-mia, notes-best.tex-rgmia}, by making use of the convolution theorem for Laplace transform and the logarithmically convex properties of the function $q_{\alpha,\beta}(x)$ on $(0,\infty)$, an alternative proof of Theorem~\ref{egp-thm-mon} was supplied.

\subsubsection{Qi-Guo's proof}
In \cite{notes-best-new-proof.tex}, by considering monotonic properties of the function
\begin{equation}\label{Qs,t;lambda(u)}
Q_{s,t;\lambda}(u)=q_{\alpha,\beta}(u)q_{\alpha,\beta}(\lambda-u),\quad \lambda\in\mathbb{R}
\end{equation}
and still employing the convolution theorem for Laplace transform, Theorem~\ref{egp-thm-mon} was completely verified again.

\begin{rem}
For more information on the function $q_{\alpha,\beta}(t)$ and its applications, please refer to \cite{bounds-two-gammas.tex, mon-element-exp.tex-rgmia, comp-mon-element-exp.tex, notes-best-new-proof.tex, mon-element-exp-final.tex, sandor-gamma-3-note.tex-final} and related references therein.
\end{rem}

\subsubsection{Qi's proof}\label{sec-6.1.2}

In \cite{notes-best-simple-open.tex, notes-best-simple.tex}, the complete monotonic properties of the function in the right-hand side of~\eqref{z-2-deriv} were established as follows.

\begin{thm}\label{byproduct}
Let $s$ and $t$ be two real numbers and $\alpha=\min\{s,t\}$. Define
\begin{equation}\label{Delta-dfn}
\Delta_{s,t}(x)=\begin{cases}\bigg[\dfrac{\psi(x+t) -\psi(x+s)}{t-s}\bigg]^2
+\dfrac{\psi'(x+t)-\psi'(x+s)}{t-s},&s\ne t\\
[\psi'(x+s)]^2+\psi''(x+s),&s=t
\end{cases}
\end{equation}
on $x\in(-\alpha,\infty)$. Then the functions $\Delta_{s,t}(x)$ for
$|t-s|<1$ and $-\Delta_{s,t}(x)$ for $|t-s|>1$ are completely
monotonic on $x\in(-\alpha,\infty)$.
\end{thm}

Since the complete monotonicity of the functions $\Delta_{s,t}(x)$ and $-\Delta_{s,t}(x)$ mean the positivity and negativity of the function $\Delta_{s,t}(x)$, an alternative proof of Theorem~\ref{egp-thm-mon} was provided once again.
\par
One of the key tools or ideas used in the proofs of Theorem~\ref{byproduct} is the following simple but specially successful conclusion: If $f(x)$ is a function defined on an infinite interval $I\subseteq\mathbb{R}$ and satisfies $\lim_{x\to\infty}f(x)=\delta$ and $f(x)-f(x+\varepsilon)>0$ for $x\in I$ and some fixed number $\varepsilon>0$, then $f(x)>\delta$ on $I$.
\par
It is clear that Theorem~\ref{byproduct} is a generalization of the inequality~\eqref{positivity}.

\subsection{Complete monotonicity of divided differences}
In order to prove the above Theorem~\ref{byproduct}, the following complete monotonic properties of a function related to a divided difference of the psi function were discovered in \cite{notes-best-simple.tex}.

\begin{thm}\label{byproduct-1}
Let $s$ and $t$ be two real numbers and $\alpha=\min\{s,t\}$. Define
\begin{equation}\label{differen-ineq}
\delta_{s,t}(x)=
\begin{cases}
\dfrac{\psi(x+t)-\psi(x+s)}{t-s}-\dfrac{2x+s+t+1}{2(x+s)(x+t)},&s\ne t\\[0.8em]
\psi'(x+s)-\dfrac1{x+s}-\dfrac1{2(x+s)^2},&s=t
\end{cases}
\end{equation}
on $x\in(-\alpha,\infty)$. Then the functions $\delta_{s,t}(x)$ for
$|t-s|<1$ and $-\delta_{s,t}(x)$ for $|t-s|>1$ are completely
monotonic on $x\in(-\alpha,\infty)$.
\end{thm}

To the best of our knowledge, the complete monotonicity of functions involving divided differences of the psi and polygamma functions were investigated first in~\cite{notes-best-simple-equiv.tex-RGMIA, notes-best-simple-equiv.tex, notes-best-simple-open.tex, notes-best-simple.tex}.

\subsection{Inequalities for sums}

As consequences of proving Theorem~\ref{byproduct-1} along a different approach from \cite{notes-best-simple.tex}, the following algebraic inequalities for sums were procured in \cite{notes-best-simple-equiv.tex-RGMIA, notes-best-simple-equiv.tex} accidentally.

\begin{thm}\label{sum-ineq-thm}
Let $k$ be a nonnegative integer and $\theta>0$ a constant.
\par
If $a>0$ and $b>0$, then
\begin{equation}\label{25inequal}
\sum_{i=0}^k\frac{1}{({a}+\theta)^{i+1} ({b}+\theta)^{k-i+1}}
+\sum_{i=0}^k\frac{1}{{a}^{i+1}{b}^{k-i+1}}
>2\sum_{i=0}^k\frac{1}{({a}+\theta)^{i+1}{b}^{k-i+1}}
\end{equation}
holds for $b-a>-\theta$ and reveres for $b-a<-\theta$.
\par
If $a<-\theta$ and $b<-\theta$, then inequalities
\begin{equation}\label{26inequal}
\sum_{i=0}^{2k}\frac{1}{({a}+\theta)^{i+1} ({b}+\theta)^{2k-i+1}}
+\sum_{i=0}^{2k}\frac{1}{{a}^{i+1}{b}^{2k-i+1}}
>2\sum_{i=0}^{2k}\frac{1}{({a}+\theta)^{i+1}{b}^{2k-i+1}}
\end{equation}
and
\begin{equation}\label{27inequal}
\sum_{i=0}^{2k+1}\frac{1}{({a}+\theta)^{i+1} ({b}+\theta)^{2k-i+2}}
+\sum_{i=0}^{2k+1}\frac{1}{{a}^{i+1}{b}^{2k-i+2}}
<2\sum_{i=0}^{2k+1}\frac{1}{({a}+\theta)^{i+1}{b}^{2k-i+2}}
\end{equation}
hold for $b-a>-\theta$ and reverse for $b-a<-\theta$.
\par
If $-\theta<a<0$ and $-\theta<b<0$, then inequality~\eqref{26inequal} holds
and inequality~\eqref{27inequal} is valid for $a+b+\theta>0$ and is reversed
for $a+b+\theta<0$.
\par
If $a<-\theta$ and $b>0$, then inequality~\eqref{26inequal} holds and
inequality~\eqref{27inequal} is valid for $a+b+\theta>0$ and is reversed for
$a+b+\theta<0$.
\par
If $a>0$ and $b<-\theta$, then inequality~\eqref{26inequal} is reversed and
inequality~\eqref{27inequal} holds for $a+b+\theta<0$ and reverses for
$a+b+\theta>0$.
\par
If $b=a-\theta$, then inequalities~\eqref{25inequal}, \eqref{26inequal} and
\eqref{27inequal} become equalities.
\end{thm}

Moreover, the following equivalent relation between the inequality~\eqref{25inequal} and Theorem~\ref{byproduct-1} was found in \cite{notes-best-simple-equiv.tex-RGMIA, notes-best-simple-equiv.tex}.

\begin{thm}\label{equiv-thm}
The inequality~\eqref{25inequal} for positive numbers $a$ and $b$ is equivalent to Theorem~\ref{byproduct-1}.
\end{thm}

\subsection{Recent advances}
Recently, some applications, extensions and generalizations of the above Theorem~\ref{byproduct}, Theorem~\ref{byproduct-1}, Theorem~\ref{sum-ineq-thm} and related conclusions have been investigated in several coming published manuscripts such as \cite{simple-equiv-simple-rev.tex, AAM-Qi-09-PolyGamma.tex}. For example, Theorem~\ref{CMDT-divided-thm} stated in Remark~\ref{kazarinoff-rem-5} has been obtained in \cite{Comp-Mon-Digamma-Trigamma-Divided.tex}.

\subsection*{Acknowledgements}
This article was ever reported on Thursday 24 July 2008 as a talk in the seminar held at the RGMIA, School of Computer Science and Mathematics, Victoria University, Australia, while the author was visiting the RGMIA between March 2008 and February 2009 by the grant from the China Scholarship Council. The author would like to express many thanks to Professors Pietro Cerone and Server S.~Dragomir and other local colleagues for their invitation and hospitality throughout this period.


\begin{thebibliography}{99}

\bibitem{abram}
M. Abramowitz and I. A. Stegun (Eds), \textit{Handbook of Mathematical Functions with Formulas, Graphs, and Mathematical Tables}, National Bureau of Standards, Applied Mathematics Series \textbf{55}, 9th printing, Washington, 1970.

\bibitem{Alzer-Math-Nachr-2001}
H. Alzer, \textit{Sharp bounds for the ratio of $q$-gamma functions}, Math. Nachr. \textbf{222} (2001), no.~1, 5\nobreakdash--14.

\bibitem{Atanassov}
R. D. Atanassov and U. V. Tsoukrovski, \textit{Some properties of a class of logarithmically completely monotonic functions}, C. R. Acad. Bulgare Sci. \textbf{41} (1988), no.~2, 21\nobreakdash--23.

\bibitem{CBerg}
C. Berg, \textit{Integral representation of some functions related to the gamma function}, Mediterr. J. Math. \textbf{1} (2004), no.~4, 433\nobreakdash--439.

\bibitem{201-05-JIPAM}
Ch.-P. Chen, \textit{Monotonicity and convexity for the gamma function}, J. Inequal. Pure Appl. Math. \textbf{6} (2005), no.~4, Art.~100; Available online at \url{http://jipam.vu.edu.au/article.php?sid=574}.

\bibitem{Chu-monthly}
J. T. Chu, \textit{A modified Wallis product and some applications}, Amer. Math. Monthly \textbf{69} (1962), no.~5, 402\nobreakdash--404.

\bibitem{dutka}
J. Dutka, \textit{On some gamma function inequalities}, SIAM J. Math. Anal. \textbf{16} (1985), 180\nobreakdash--185.

\bibitem{GURLAND}
J. Gurland, \textit{On Wallis' formula}, Amer. Math. Monthly \textbf{63} (1956), 643\nobreakdash--645.

\bibitem{egp}
N. Elezovi\'c, C. Giordano and J. Pe\v{c}ari\'c, \textit{The best bounds in Gautschi's inequality}, Math. Inequal. Appl. \textbf{3} (2000), 239\nobreakdash--252.

\bibitem{gaut}
W. Gautschi, \textit{Some elementary inequalities relating to the gamma and incomplete gamma function}, J. Math. Phys. \textbf{38} (1959/60), 77\nobreakdash--81.

\bibitem{Giordano-Laforgia-01-jcam}
C. Giordano and A. Laforgia, \textit{Inequalities and monotonicity properties for the gamma function}, J. Comput. Appl. Math. \textbf{133} (2001) 387\nobreakdash--396.

\bibitem{GLP-98}
C. Giordano, A. Laforgia and J. Pe\v{c}ari\'c, \textit{Monotonicity properties for some functions involving the ratio of two gamma functions}, in: A. Bellacicco, A. Laforgia (Eds.), Funzioni Speciali e Applicazioni, Franco Angeli, Milano, 1998, 35\nobreakdash--42.

\bibitem{unify-kershaw}
C. Giordano, A. Laforgia and J. Pe\v{c}ari\'c, \textit{Unified treatment of Gautschi-Kershaw type inequalities for the gamma function}, Proceedings of the VIIIth Symposium on Orthogonal Polynomials and Their Applications (Seville, 1997). J. Comput. Appl. Math. \textbf{99} (1998), no.~1-2, 167\nobreakdash--175.

\bibitem{Kazarinoff-56}
D. K. Kazarinoff, \textit{On Wallis' formula}, Edinburgh Math. Notes \textbf{1956} (1956), no.~40, 19\nobreakdash--21.

\bibitem{kershaw}
D. Kershaw, \textit{Some extensions of W. Gautschi's inequalities for the gamma function}, Math. Comp. \textbf{41} (1983), 607\nobreakdash--611.

\bibitem{laforgia-mc-1984}
A. Laforgia, \textit{Further inequalities for the gamma function}, Math. Comp. \textbf{42} (1984), no.~166, 597\nobreakdash--600.

\bibitem{Lazarevic}
I. Lazarevi\'c and A. Lupa\c{s}, \textit{Functional equations for Wallis and Gamma functions}, Publ. Elektrotehn. Fak. Univ. Beograd. Ser. Electron. Telecommun. Automat. No.~\textbf{461-497} (1974), 245\nobreakdash--251.

\bibitem{lorch-ultra}
L. Lorch, \textit{Inequalities for ultraspherical polynomials and the gamma function}, J. Approx. Theory \textbf{40} (1984), no.~2, 115\nobreakdash--120.

\bibitem{Merkle-JMAA-99}
M. Merkle, \textit{Representations of error terms in Jensen's and some related inequalities with applications}, J. Math. Anal. Appl. \textbf{231} (1999), 76\nobreakdash--90.

\bibitem{mit}
D. S. Mitrinovi\'c, \textit{Analytic Inequalities}, Springer-Verlag, 1970.

\bibitem{mpf-1993}
D. S. Mitrinovi\'c,  J. E. Pe\v{c}ari\'c and A. M. Fink, \textit{Classical and New Inequalities in Analysis}, Kluwer Academic Publishers, 1993.

\bibitem{Palumbo-98-jcam}
B. Palumbo, \textit{A generalization of some inequalities for the gamma function}, J. Comput. Appl. Math. \textbf{88} (1998), no.~2, 255\nobreakdash--268.

\bibitem{notes-best-simple-equiv.tex-RGMIA}
F. Qi, \textit{A completely monotonic function involving divided difference of psi function and an equivalent inequality involving sum}, RGMIA Res. Rep. Coll. \textbf{9} (2006), no.~4, Art.~5; Available online at \url{http://www.staff.vu.edu.au/rgmia/v9n4.asp}.

\bibitem{notes-best-simple-equiv.tex}
F. Qi, \textit{A completely monotonic function involving the divided difference of the psi function and an equivalent inequality involving sums}, ANZIAM J. \textbf{48} (2007), no.~4, 523\nobreakdash--532.

\bibitem{notes-best-simple-open.tex}
F. Qi, \textit{A completely monotonic function involving divided differences of psi and polygamma functions and an application}, RGMIA Res. Rep. Coll. \textbf{9} (2006), no.~4, Art.~8; Available online at \url{http://www.staff.vu.edu.au/rgmia/v9n4.asp}.

\bibitem{bounds-two-gammas.tex}
F. Qi, \textit{Bounds for the ratio of two gamma functions}, RGMIA Res. Rep. Coll. \textbf{11} (2008), no.~3, Art.~1; Available online at \url{http://www.staff.vu.edu.au/rgmia/v11n3.asp}.

\bibitem{mon-element-exp.tex-rgmia}
F. Qi, \textit{Monotonicity and logarithmic convexity for a class of elementary functions involving the exponential function}, RGMIA Res. Rep. Coll. \textbf{9} (2006), no.~3, Art.~3; Available online at \url{http://www.staff.vu.edu.au/rgmia/v9n3.asp}.

\bibitem{Comp-Mon-Digamma-Trigamma-Divided.tex}
F. Qi, \textit{Necessary and sufficient conditions for a function involving divided differences of the di- and tri-gamma functions to be completely monotonic}, submitted.

\bibitem{notes-best-simple.tex}
F. Qi, \textit{The best bounds in Kershaw's inequality and two completely monotonic functions}, RGMIA Res. Rep. Coll. \textbf{9} (2006), no.~4, Art.~2; Available online at \url{http://www.staff.vu.edu.au/rgmia/v9n4.asp}.

\bibitem{comp-mon-element-exp.tex}
F. Qi, \textit{Three-log-convexity for a class of elementary functions involving exponential function}, J. Math. Anal. Approx. Theory \textbf{1} (2006), no.~2, 100\nobreakdash--103.

\bibitem{simple-equiv-simple-rev.tex}
F. Qi, P. Cerone and S. S. Dragomir, \textit{Complete monotonicity results of divided difference of psi functions and new bounds for ratio of two gamma functions}, submitted.

\bibitem{notes-best-new-proof.tex}
F. Qi and B.-N. Guo, \textit{An alternative proof of Elezovi\'c-Giordano-Pe\v{c}ari\'c's theorem}, submitted.

\bibitem{minus-one}
F. Qi and B.-N. Guo, \textit{Complete monotonicities of functions involving the gamma and digamma functions}, RGMIA Res. Rep. Coll. \textbf{7} (2004), no.~1, Art.~8, 63\nobreakdash--72; Available online at \url{http://www.staff.vu.edu.au/rgmia/v7n1.asp}.

\bibitem{AAM-Qi-09-PolyGamma.tex}
F. Qi and B.-N. Guo, \textit{Necessary and sufficient conditions for functions involving the tri- and tetra-gamma functions to be completely monotonic}, Adv. Appl. Math. (2009), in press.

\bibitem{Infinite-family-Digamma.tex}
F. Qi and B.-N. Guo, \textit{Sharp inequalities for the psi function and harmonic numbers}, submitted.

\bibitem{mon-element-exp-final.tex}
F. Qi and B.-N. Guo, \textit{Properties and applications of a function involving exponential functions}, Commun. Pure Appl. Anal. (2009), in press.

\bibitem{sandor-gamma-3-note.tex-final}
F. Qi and B.-N. Guo, \textit{Wendel's and Gautschi's inequalities: Refinements, extensions, and a class of logarithmically completely monotonic functions}, Appl. Math. Comput. \textbf{205} (2008), no.~1, 283\nobreakdash--292; Available online at \url{http://dx.doi.org/10.1016/j.amc.2008.07.005}.

\bibitem{sandor-gamma-3-note.tex}
F. Qi and B.-N. Guo, \textit{Wendel-Gautschi-Kershaw's inequalities and sufficient and necessary conditions that a class of functions involving ratio of gamma functions are logarithmically completely monotonic}, RGMIA Res. Rep. Coll. \textbf{10} (2007), no.~1, Art.~2; Available online at \url{http://www.staff.vu.edu.au/rgmia/v10n1.asp}.

\bibitem{notes-best.tex-mia}
F. Qi, B.-N. Guo and Ch.-P. Chen, \textit{The best bounds in Gautschi-Kershaw inequalities}, Math. Inequal. Appl. \textbf{9} (2006), no.~3, 427\nobreakdash--436.

\bibitem{notes-best.tex-rgmia}
F. Qi, B.-N. Guo and Ch.-P. Chen, \textit{The best bounds in Gautschi-Kershaw inequalities}, RGMIA Res. Rep. Coll. \textbf{8} (2005), no.~2, Art.~17; Available online at \url{http://www.staff.vu.edu.au/rgmia/v8n2.asp}.

\bibitem{Rao-Uppuluri}
V. R. Rao Uppuluri, \textit{On a stronger version of Wallis' formula}, Pacific J. Math. \textbf{19} (1966), no.~1, 183\nobreakdash--187.

\bibitem{waston}
G. N. Watson, \textit{A note on gamma functions}, Proc. Edinburgh Math. Soc. \textbf{11} (1958/1959), no.~2, Edinburgh Math Notes No.~42 (misprinted 41) (1959), 7\nobreakdash--9.

\bibitem{LaplaceTransform.html}
E. W. Weisstein, \textit{Laplace Transform}, From MathWorld---A Wolfram Web Resource; Available online at \url{http://mathworld.wolfram.com/LaplaceTransform.html}.

\bibitem{WallisFormula.html}
E. W. Weisstein, \textit{Wallis Cosine Formula}, From MathWorld\nobreakdash---A Wolfram Web Resource; Available online at \url{http://mathworld.wolfram.com/WallisFormula.html}.

\bibitem{wendel}
J. G. Wendel, \textit{Note on the gamma function}, Amer. Math. Monthly \textbf{55} (1948), no.~9, 563\nobreakdash--564.

\bibitem{widder}
D. V. Widder, \textit{The Laplace Transform}, Princeton University Press, Princeton, 1946.

\end{thebibliography}
\end{document}